\documentclass{amsart}
\usepackage{amscd,amsmath}

\newcommand{\ov}{\overline}
\newcommand{\Ker}{\mathop\mathrm{Ker\,}}

\newcommand{\dif}{\mathrm{d}}
\newcommand{\tr}{\mathrm{trace}}
\newcommand{\CC}{\mathbb{C}}

\newcommand{\HH}{\mathcal{H}}
\newcommand{\VV}{\mathcal{V}}
\newcommand{\FF}{\mathcal{F}} 
\newcommand{\DD}{\mathcal{D}} 
\newcommand{\RR}{\mathbb{R}}

\numberwithin{equation}{section}

\newtheorem{te}{Theorem}[section]
\newtheorem{pr}{Proposition}[section]
\newtheorem{co}{Corollary}[section]
\newtheorem{lm}{Lemma}[section]

\theoremstyle{definition}
\newtheorem{de}{Definition}[section]    
\newtheorem{re}{Remark}[section]
\newtheorem{ex}{Example}[section]

\begin{document}  

\title{Pseudo-harmonic morphisms with low dimensional fibers}
\author{Radu Slobodeanu}
\address{Faculty of Physics, Bucharest University - 405 Atomistilor str., CP Mg-11, RO-76900 Bucharest, Romania.}
\email{slobyr@yahoo.com}

\date{December 2007 for this revisited version of the paper published in {\it Rend. Circ. Mat. Palermo}, Serie II, Tomo \textbf{LV} (2006), 5 - 20.}

\maketitle

\begin{abstract}
We characterize pseudo-harmonic morphisms from a Riemannian manifold to a Hermitian manifold as pseudo horizontally weakly conformal maps with an additional property. We study to what extent we can (locally) describe these submersive pseudo-harmonic morphisms \textit{via} the foliation given by the kernel of the associated $f$-structure.

In a second part, we point out that, in the case of pseudo-harmonic morphisms with one and two-dimensional fibers, the induced $f$-structure gives rise to an almost contact, respectively almost complex structure on the domain. We give criteria for normality and integrability of these structures and we show how these two particular cases are interrelated. 
\end{abstract}

\section{Introduction}

Harmonic morphisms are those maps between Riemannian manifolds that pull local harmonic functions on the codomain back to local harmonic functions on the domain, \cite{ud}. Analogously, one can define \textit{pseudo-harmonic morphisms} (PHM) as maps from a Riemannian manifold to an almost complex one, that pull local \textit{holomorphic} functions on the codomain back to local harmonic functions on the domain.

If the codomain is K\"ahler, pseudo-harmonic morphisms have a nice description similar to Fuglede-Ishihara's characterization for  harmonic morphisms, cf. \cite {lub}: a map is a pseudo-harmonic morphism if and only if it is both pseudo horizontally weakly conformal and harmonic. Recall that a map $\varphi: (M, g) \longrightarrow (N, J, h)$ from a Riemannian manifold to an almost Hermitian one is \emph{pseudo horizontally weakly conformal} (PHWC) if, by definition, \cite {cen}, \cite {lub}:

\begin{equation} \label{PHWC}
[d\varphi \circ d\varphi^{*}, J]=0
\end{equation}

\noindent where $d\varphi^{*}$ stands for the \textit{adjoint map}, $d\varphi^{*}_x: T_{\varphi(x)}N \rightarrow T_x M$, characterized by: $g(X, d\varphi^{*}_x (E)) = h(d\varphi_x(X), E)$, $\forall X \in T_x M, E \in T_{\varphi(x)}N$. Note that horizontally weakly conformal (HWC) maps are obviously PHWC, as  in their case, $d\varphi \circ d\varphi^{*} = \lambda^{2} id$. In particular, the harmonic morphisms onto a K\"ahler manifold are also pseudo harmonic morphisms.

If the PHWC condition is satisfied, according to \cite[Prop.7]{lub}, a metric $f$-structure $F^{\varphi}$ is induced on $(M,g)$ such that $\varphi$ becomes $(F^{\varphi}, J)$-holomorphic. The geometric meaning of \eqref{PHWC} becomes more transparent if we notice that, when $\varphi$ is submersive, its differential induces an almost complex structure on the horizontal bundle, by taking $J_{\mathcal{H}} = d \varphi \vert _{\mathcal{H}}^{-1} \circ J \circ d\varphi \vert_{\mathcal{H}}$. Then one can prove that the PHWC condition is equivalent to the compatibility of $J_{\mathcal{H}}$ with the domain metric $g$ (in this case $\varphi$ becomes \textit{horizontally holomorphic}, according to the terminology used in \cite{panti}). 

A classical example of PHWC map is any stable harmonic map to an irreducible Hermitian symmetric space of compact type, cf. \cite {bur}.
 
The analogue of horizontal homothety in this context was introduced by M.A. Aprodu, M. Aprodu and V. Br\^{\i}nz\u anescu in \cite {aab}. A \emph{pseudo horizontally homothetic} (PHH) map is a PHWC map $\varphi$ which satisfies:

\begin{equation} \label{PHH}
[d\varphi \circ \nabla^{M}_{X} \circ d\varphi^{*}, J]=0, \forall X\in\Gamma(\mathcal{H})
\end{equation}

In turn, this condition means that $J_{\mathcal{H}}$ is parallel (with respect to $\nabla^{\mathcal{H}}$) in horizontal directions, so satisfies a \emph{transversal K\"ahler} condition. Any PHH harmonic submersion onto a K\"ahler manifold exhibits a particularly nice geometric property: it pulls complex submanifolds back to minimal submanifolds, cf. \cite {aab}. Moreover, in \cite {ma}, it is shown that PHH harmonic submersions (from a compact manifold) are (weakly) stable (as holomorphic maps between K\"ahler manifolds do). Further properties and examples of PHH harmonic submersions can be found in \cite {aa}, \cite{brz}.

\medskip
This paper is organized as follows. In the next section we study some properties of PHWC submersions to a Hermitian manifold and we describe the foliations given locally by the fibers of this kind of submersions. Section 3 is devoted to pseudo harmonic morphisms in the general sense mentioned above. We check if the canonical foliation of an (almost) contact manifold can be given locally by the fibers of a PHM. In section 4 and 5 we suppose that the fibers of a submersive PHM are 1-dim., 2-dim. respectively. We point out that the induced $f$-structure gives rise to an almost contact, respectively an almost complex structure on the domain. We give criteria for normality and integrability in each case and we stress the interrelation between these two particular cases. The paper ends with an Appendix and few comments gathering both the motivation and conclusion of this work.

Note that our approach differs from the one proposed in \cite{lub}, where the codomain is always supposed K\"ahler in order to speak about PHM's (besides, PHM's are defined by the conditions of being harmonic and PHWC). In what follows, we shall stress the difference by calling a map \textit{strong PHM} when it is PHM with K\"ahler target, as in \cite{lub}.

\section{PHWC submersions and associated foliations} 

\subsection{PHWC submersions} As it has been shown in \cite{lub}, any PHWC map $\varphi: (M, g) \longrightarrow (N, J, h)$ from a Riemannian manifold to a Hermitian one induces a metric $f$-structure, $F$, on $M$, with respect to which $\varphi$ becomes $(F, J)$-holomorphic. When $\varphi$ is submersive, with $\HH$ and $\VV$ horizontal and vertical distributions respectively, $F$ extends the induced almost complex structure on $\HH$: $F | _{\mathcal{H}}=J_{\mathcal{H}}=d\varphi\vert _{\mathcal{H}}^{-1} \circ J \circ d\varphi \vert_{\mathcal{H}}$ and $F | _{\mathcal{V}}=0$. Then the PHWC condition is equivalent to the compatibility of $J_{\HH}$ with the domain metric $g$ (i.e. $F$ is a metric $f$-structure). 
Indeed, starting with the remark that $d\varphi^*(E)$ is horizontal for all $E$ tangent to $N$, we have
\begin{equation*}
\begin{split}
g(J_{\mathcal{H}}X, d\varphi^* E)&=h(d\varphi(J_{\HH}X), E)=
h(J (d\varphi X), E)=-h(d\varphi X, JE)=-g(X, d\varphi^* JE)\\
&=-g(X, d\varphi^{-1}(d\varphi(d\varphi^* JE)))=
-g(X, d\varphi^{-1}(J(d\varphi(d\varphi^* E)))\\
&=-g(X, J_{\mathcal{H}}d\varphi^* E)
\end{split}
\end{equation*}
and because $d\varphi^*:TN \rightarrow \HH$ is an isomorphism, the compatibility statement follows.
\begin{lm}
Let F be the f-structure on M, naturally induced by the PHWC submersion 
$\varphi: (M, g) \longrightarrow (N, J, h)$. 
Then the following relation holds:
\begin{equation}
((\mathcal{L}_{V} F)X)^{\HH} = 0, \ \forall X \in \Gamma(\HH), V \in \Gamma(\VV),
\end{equation}
where  $\HH$ and $\VV$ are, respectively, the horizontal and vertical distributions of the submersion $\varphi$.
\end{lm}

\begin{proof}
By hypothesis, we have: $d\varphi \vert_{\mathcal{H}} \circ J_{\mathcal{H}} = J \circ d\varphi \vert_{\mathcal{H}}$.
It is easy to see that if $X \in \Gamma(\mathcal{H})$ is projectable (i.e. $\varphi$-related with a vector field, $\check{X}$, on $N$), then $FX=J_{\mathcal{H}}X$ is projectable too (i.e. $\varphi$-related with  $J\check{X}$, on $N$). Therefore $[V,X]$, $[V,FX]$ are both vertical, so $\left((\mathcal{L}_{V}F)X \right)^{\HH}=([V,FX]-F[V,X])^{\HH}=0$ is true for any projectable vector field $X$, so it is true in general (we always have a local frame of projectable vector fields for the horizontal distribution of a submersion).
\end{proof}

\subsection{PHWC foliations}
First, let us recall some basic notions in the study of distributions / foliations on a Riemannian manifold $(M^{m}, g)$. Let $\mathcal{V}$ and $\mathcal{H}$ be two orthogonal complementary distributions of dimension $m-n$ and $n$, respectively. The exponent $\mathcal{V}$ or $\mathcal{H}$ will indicate the orthogonal projections onto these distributions. The second fundamental form $B^{\mathcal{H}}$ and the integrability tensor $I^{\mathcal{H}}$ of $\mathcal{H}$ are defined by:
$$
B^{\mathcal{H}}(X, Y)=\frac{1}{2}\left(\nabla_{X}Y+\nabla_{Y}X\right)^{\mathcal{V}},\quad
I^{\mathcal{H}}(X, Y)=[X, Y]^{\mathcal{V}}, \quad \forall X, Y \in \Gamma(\mathcal{H}).
$$
Whether $\mathcal{H}$ is integrable or not, {\it the mean curvature} of $\mathcal{H}$ is, by definition:
$$
\mu^{\mathcal{H}}=\frac{1}{n}\text{trace}B^{\mathcal{H}}.
$$
The \emph{Bott partial connection} 
$\stackrel{\circ}{\nabla}=\stackrel{\circ}{\nabla^{\mathcal{H}}}$
on $\mathcal{H}$ is the map $\stackrel{\circ}{\nabla}: \Gamma(\mathcal{V}) \times
\Gamma(\mathcal{H}) \longrightarrow \Gamma(\mathcal{H})$ defined by
$$
\stackrel{\circ}{\nabla}_{V}X=[V, X]^{\mathcal{H}}.
$$
Extending Bott connections to horizontal tensor fields (i.e. sections of $\otimes^r \HH \otimes \otimes^s \HH^*$), we say that $\sigma$ is \textit{basic} if $\stackrel{\circ}{\nabla}_{V} \sigma=0$.

\begin{de} 
Let $(M^{m}, g)$ a Riemannian manifold. A distribution $\mathcal{V}$ on $M$ is called \emph{transversally (almost) Hermitian} if the complementary distribution $\mathcal{H}$ admits an (almost) complex structure $J_{\mathcal{H}}$ compatible with the metric $g$. In particular, the codimension of $\VV$ must be even.
\end{de}

\begin{pr}
Let $(M^{m}, g)$ a Riemannian manifold equipped with a transversally almost Hermitian distribution $\mathcal{V}$. Then the following relations hold good:

\medskip
$(i)$ $g\left( (\nabla_{E}J_{\mathcal{H}})X, X \right)=0, \ \forall E\in\Gamma(TM), X\in\Gamma(\mathcal{H})$

\medskip
$(ii)$ $g\left( (\nabla_{E}J_{\mathcal{H}})X, J_{\mathcal{H}}X \right)=0, \ \forall E\in\Gamma(TM), X\in\Gamma(\mathcal{H})$

\medskip
$(iii)$ $g\left( (\stackrel{\circ}{\nabla}_{V}J_{\mathcal{H}})X, X \right)=
-(\mathcal{L}_{V}g)(X, J_{\mathcal{H}}X), \ \forall V\in\Gamma(\mathcal{V}), 
X\in\Gamma(\mathcal{H})$

\medskip
$(iv)$ $g\left( (\stackrel{\circ}{\nabla}_{V}J_{\mathcal{H}})X, J_{\mathcal{H}}X \right)=
-\frac{1}{2}(\mathcal{L}_{V}g)(J_{\mathcal{H}}X, J_{\mathcal{H}}X)
+\frac{1}{2}(\mathcal{L}_{V}g)(X, X), \ \forall V\in\Gamma(\mathcal{V})$, 
$X\in\Gamma(\mathcal{H})$.
\end{pr}

\begin{proof}
The proof is essentially the same as for \cite[Prop. 2.5.16] {ud}.
\end{proof}

Now let $\VV$ be integrable so that it defines a foliation $\FF$. If $\VV$ is transversally almost Hermitian in the sense of the above definition, then $\FF$ is a {\it transversally almost Hermitian foliation} if, in addition, the almost complex structure $J_{\HH}$ and the metric $g$ are basic tensors. 

We have seen that $J_{\HH}$ induced by a PHWC submersion is basic (i.e. $\stackrel{\circ}{\nabla}_{V}J_{\HH}=0$). However, the metric is not necessarily basic, so in general, PHWC submersions does not induce transversally almost Hermitian foliations on the domain. 

\begin{de}
A {\it PHWC foliation} on a Riemannian manifold, is a foliation for which the normal bundle has a basic almost complex structure compatible with the ambient metric. 
\end{de}

Therefore transversally almost Hermitian foliations are Riemannian PHWC foliations. Note also that any PHWC foliation $\FF$ on $(M,g)$ induces a metric $f$-structure on $M$ by $F \vert _{\HH}=J_{\HH}$, $F \vert _{\mathcal{V}}=0$ (we shall call it the \textit{induced $f$-structure}), where we denote by $\VV$ the tangent distribution to $\FF$  and by $\HH$ the normal bundle identified with $\VV^{\perp}$ due to the metric on the ambient space. For further details on the foliations we refer to \cite{ton}.

As a direct consequence of Proposition 2.1 we point out that:
\begin{co}
For a PHWC foliation $\FF$ we have the following equivalent relations:

\medskip
$(i)$ $(\mathcal{L}_{V}g)(J_{\mathcal{H}}X, J_{\mathcal{H}}Y)
=(\mathcal{L}_{V}g)(X, Y), \ \forall V\in\Gamma(\mathcal{V}), X, Y \in\Gamma(\mathcal{H})$. 

\medskip
$(ii)$ $B^{\mathcal{H}}(J_{\mathcal{H}}X, J_{\mathcal{H}}Y)=
B^{\mathcal{H}}(X, Y), \ \forall X, Y \in\Gamma(\mathcal{H})$. 
\end{co}

\begin{proof}
The equivalence is assured by the following relation:
$$
(\mathcal{L}_{V}g)(X, Y)=
-2g(B^{\mathcal{H}}(X, Y), V), \quad \forall V \in \Gamma(\mathcal{V}), X, Y \in\Gamma(\mathcal{H}).
$$
\end{proof}

\begin{pr} \emph{(Characterization of PHWC foliations)}
A codimension 2n foliation $\mathcal{F}$ on a Riemannian manifold is a PHWC foliation if and only if for each $\mathcal{F}$-simple open set U, the leaf space $U / \mathcal{F}$ can be endowed with an almost Hermitian structure with respect to which the natural projection
$U \longrightarrow U/ \mathcal{F}$ is a PHWC submersion.
\end{pr}

\begin{proof}
The "only if" part has been already remarked, so we have to prove only the "if" part.

Let $\FF$ be a PHWC foliation and $U$ a $\mathcal{F}$-simple open set in $M$. Denote $U / \FF$ by $\ov{U}$. Because the almost complex structure that we have on $\HH$ is basic, one can construct an almost complex structure $J$ on $\ov{U}$ by:
$J=d\varphi \circ J_{\HH} \circ d\varphi^{-1}$, where $\varphi$ denotes the standard projection and $d\varphi$ is restricted to $\HH$. It is obvious that $J^{2}=-id$ and that we can take a Riemannian metric $h$ on $\ov{U}$, compatible with $J$ (so it exists an almost Hermitian structure on $\ov{U}$).

Now we shall verify the PHWC condition in this context, that is:

$d\varphi \circ d\varphi^{*} \circ J=J \circ d\varphi \circ d\varphi^{*} \ 
\Leftrightarrow
\ d\varphi^{*} \circ J \circ (d\varphi^{*})^{-1} = d\varphi^{-1} \circ J \circ d\varphi$

$\Leftrightarrow \ d\varphi^{*} \circ J \circ (d\varphi^{*})^{-1}=J_{\HH}$.

In terms of associated matrices, with obvious notations, this translates into:
$(G^{-1}\Phi^{t}H)J(H^{-1}(\Phi^{t})^{-1}G)=J_{\mathcal{H}}$.
But the compatibility of $J$ with the metric $h$ means simply: $HJ=JH$, so the above relation becomes:
$\left(\Phi^{-1}J^{t}\Phi\right)^{t}G=GJ_{\mathcal{H}} \Leftrightarrow  J_{\mathcal{H}}G = GJ_{\mathcal{H}}$,
which is clearly satisfied by hypothesis (compatibility of $J_{\mathcal{H}}$ with $g$).
\end{proof}

\bigskip

Let us notice that, generally, the induced almost Hermitian structure in the above proof is neither integrable, nor K\"ahler (as requested for strong pseudo-harmonic morphisms).

\subsection{Integrability of the induced $f$-structure}
Note firstly an important feature of PHWC submersions, that a certain amount of integrability of $J$ is inherited by the induced $f$-structure, $F$:

\begin{lm}  
Let $\varphi: (M,g) \longrightarrow (N^{2n}, J, h)$ be a PHWC
submersion onto a Hermitian manifold and $F$ be the f-structure induced on M.
Then: 
$$
[F, F](X, Y)^{\HH}=0, \forall X, Y \in \Gamma (\HH).
$$
Moreover, if N is K\"ahler, then on the horizontal bundle $\HH$ it is defined a nondegenerate closed F-invariant 2-form.
\end{lm}   

\begin{proof} 
By definition, $F X=d\varphi^{-1} \circ J \circ d\varphi (X), \forall X \in \Gamma(\mathcal{H})$, where $d\varphi^{-1}$ stands for the
horizontal lift application. For any basic vector fields $X, Y$
we shall have the following sequence of identities (similar to those one in \cite [Example 6.7.2]{ble}):   
\begin{equation*} 
\begin{split} 
[F, F](X, Y)& =  d\varphi^{-1}J^{2}d\varphi [X, Y]+ 
[d\varphi^{-1}Jd\varphi X, d\varphi^{-1}Jd\varphi Y]\\ 
&-d\varphi^{-1}Jd\varphi [d\varphi^{-1}Jd\varphi X, Y] 
-d\varphi^{-1}Jd\varphi [X, d\varphi^{-1}Jd\varphi Y]\\ 
&=d\varphi^{-1}J^{2}[d\varphi X, d\varphi Y]+ 
d\varphi^{-1}[Jd\varphi X, Jd\varphi Y]\\ 
&+[d\varphi^{-1}Jd\varphi X, d\varphi^{-1}Jd\varphi Y]^{\VV}\\ 
&-d\varphi^{-1}J [Jd\varphi X, d\varphi Y] 
-d\varphi^{-1}J [d\varphi X, Jd\varphi Y]\\ 
&=d\varphi^{-1}[J, J][d\varphi X, d\varphi Y]+ 
[d\varphi^{-1}Jd\varphi X, d\varphi^{-1}Jd\varphi Y]^{\VV}\\ 
&=[F X, F Y]^{\VV}, 
\end{split} 
\end{equation*} 
because $J$ is, by hypothesis, integrable. 

For the second assertion, simply take the pull-back
$\varphi^{*}\Omega$ of the K\"ahler form on $N$.
\end{proof}

We shall work with the following definition:

\begin{de}(\cite{pant})
An almost $f$-structure on $M$ is a section $F$ of $End(TM)$ such that
$F^3+F = 0$. Let $\FF = T^0 M \oplus T^{0,1} M$ where $T^0 M$ and $T^{0,1} M$ are the eigenbundles of $F$ corresponding
to 0 and -i , respectively; we say that $\FF$ is \textit{the complex distribution associated to} $F$. Then $F$ is \textit{integrable} if $\FF$ is integrable (that is, for any $X, Y \in \Gamma(\FF)$ we have $[X, Y ] \in \Gamma(\FF)$).
\end{de}

The following class of examples is based upon a well-known theorem of Ianu\c s, \cite{ia}.

\begin{ex}
If an almost contact structure is \textit{normal}, then it is \textit{integrable} according to the definition cited above. The converse is not necessarily true.
\end{ex}

We point out now that the \textit{induced $f$-structures} by a PHWC submersion onto a Hermitian manifold are integrable. Next we shall see that \textit{all} integrable $f$-structures (locally) appear in this way.

\begin{pr}
The $f$-structure induced by a PHWC submersion onto a Hermitian manifold $\varphi: (M, g) \longrightarrow (N, J, h)$ is integrable. 
\end{pr}

\begin{proof}
The proof follows directly from Lemma 2.1 and Lemma 2.2. Indeed, 

\noindent $((\mathcal{L}_{V} F)X)^{\HH} = 0$ assures us that 
$[aV, X+iFX]\in T^0 M \oplus T^{0,1} M$.
On the other hand $[F, F](X, Y)^{\HH}=0, \forall X, Y \in \Gamma (\HH)$
tells us exactly that $[X+iFX, Y+iFY]\in T^0 M \oplus T^{0,1} M$.
\end{proof}

If the fibers of $\varphi$ are 1-dim. and $M$ is orientable, we shall see that the induced $f$-structure is an almost contact metric structure, $\phi$. It is integrable but not necessarily normal. The supplementary conditions needed in order to be normal will be given in Theorem 4.1 below.

\bigskip

\begin{pr} \emph{(Characterization of integrable $f$-structures, \cite{pant})}
An almost $f$-structure $F$ on $M$ is integrable if and only if for any
$x \in M$ there exists an open neighbourhood $U$ of $x$ and a $(F,J)$-holomorphic submersion $\varphi$ from $(U, F| _U)$ onto some complex manifold $(N, J)$ such that $\Ker \mathrm{d} \varphi = T^0 M$; we say that the $f$-structure $F|_U$ is \textit{defined by} $\varphi$. A \textit{simple} $f$-structure is an $f$-structure (globally) defined by a holomorphic submersion with connected fibres.
\end{pr}

\section{Pseudo-harmonic morphisms and the foliations that produce them}

\subsection{General facts about pseudo-harmonic morphisms.} 
Recall that, \cite[Prop. 2]{lub}: 
\begin{quote}
\textit{If $(N, J, h)$ is a K\"ahler manifold, then $\varphi$ is a (strong) PHM if and only if it is PHWC and harmonic.}
\end{quote}

In this context, the condition for harmonicity of PHWC mappings is given by, \cite {mo}:

\begin{quote} 
{\it Let $\varphi: (M, g) \longrightarrow (N, J, h)$
be a PHWC map from a Riemannian to a K\"ahler manifold. Then $\varphi$ is harmonic (and therefore a strong PHM) if and only if} $F \mathrm{div} F=0$.
\end{quote}
This is simply a consequence of the particular form that the tension field of any PHWC map takes:
\begin{equation} \label{ta}
\tau(\varphi)=J \mathrm{div}^{\varphi} J - d \varphi(F \mathrm{div} F), 
\end{equation}

where $\mathrm{div}^{\varphi} J = \tr _g \varphi^{*}\nabla^N J$ and $\mathrm{div} F$ = trace $\nabla F$ is the divergence of $F$.
Locally, if we consider an adapted frame $\{e_{i}, F e_{i}, e_{\alpha}\}$ (i.e. an orthonormal frame such that  $e_{\alpha} \in \Ker{F}$), then the above relation reads: 

\begin{equation} \label{tau}
\begin{split}
\tau(\varphi)= & J\left( \sum_{i=1}^{n}(\nabla^{\varphi}_{e_{i}}J)(d\varphi(e_{i}))+ (\nabla^{\varphi}_{Fe_{i}}J)(d\varphi(Fe_{i}))\right) \\
&-d \varphi\left( \sum_{i=1}^{n}F \left[(\nabla_{e_{i}}F)(e_{i})+ (\nabla_{F e_{i}}F)(Fe_{i})\right]+(m-2n)\mu^{\mathcal{V}} \right).
\end{split} 
\end{equation}

If $F \mathrm{div} F=0$ then the $f$-structure, $F$, is called \textit{cosymplectic}. Note that this is equivalent to: $\mu^{\mathcal{V}} = -\frac{1}{(m-2n)}F\mathrm{div}^{\mathcal{H}}F$.

\medskip

Let us turn to the general setup for pseudo-harmonic morphisms, which does not involve a metric on the codomain (a general PHM is a map $\varphi: (M, g) \longrightarrow (N, J)$ from a Riemannian manifold to a complex manifold), so PHWC and harmonic properties are \textit{a priori} meaningless. Nevertheless if we consider $N$ locally endowed with the standard flat metric, the characterization of E. Loubeau \cite[Prop. 2]{lub} holds and our PHM is PHWC and harmonic (locally).

Now fix an arbitrary Hermitian metric (globally) on $(N,J)$. We have the following:

\begin{te}
Let $\varphi: (M, g) \longrightarrow (N, J, h)$ be a submersive map to a Hermitian manifold. Then $\varphi$ is a (general) pseudo-harmonic morphism to a Hermitian manifold if and only if $\varphi$ is a PHWC map with cosymplectic induced $f$-structure.
\end{te}

\begin{proof}
Suppose that $\varphi$ is PHM. We have to adapt the proof of \cite[Prop. 2]{lub}.

Let $f$ be a local holomorphic function on $N$. As $\varphi$ pulls  germs of holomorphic functions back to germs of harmonic functions, we must have $\tau(f\circ \varphi)=0$.
This translates in a local chart as:
\begin{equation}\label{comp}
\dfrac{\partial f}{\partial z ^{\gamma}}\tau^{\gamma}(\varphi)
+g^{ij}\left[\dfrac{\partial^2 f}{\partial z ^{\alpha} \partial z^{\beta}} \dfrac{\partial \varphi^{\alpha}}{\partial x ^{i}}\dfrac{\partial \varphi^{\beta}}{\partial x ^{j}} -  \,^{N}\Gamma_{\alpha \ov \beta}^{\delta}\dfrac{\partial f}{\partial z ^{\delta}}\dfrac{\partial \varphi^{\alpha}}{\partial x ^{i}}\dfrac{\partial \varphi^{\ov \beta}}{\partial x ^{j}}\right]=0.
\end{equation}

Taking $f=z^{\gamma}$ we get:
$$
\tau^{\gamma}(\varphi)-g^{ij}  \,^{N}\Gamma_{\alpha \ov \beta}^{\gamma}\dfrac{\partial f}{\partial z ^{\gamma}}\dfrac{\partial \varphi^{\alpha}}{\partial x ^{i}}\dfrac{\partial \varphi^{\ov \beta}}{\partial x ^{j}}=0, \quad \forall \gamma.
$$

Therefore \eqref{comp} reduces to: $g^{ij}\frac{\partial^2 f}{\partial z ^{\alpha} \partial z^{\beta}} \frac{\partial \varphi^{\alpha}}{\partial x ^{i}} \frac{\partial \varphi^{\beta}}{\partial x ^{j}}=0$, 
for any holomorphic function $f$. Giving particular values to $f$ we get:
$$
g^{ij} \dfrac{\partial \varphi^{\alpha}}{\partial x ^{i}}\dfrac{\partial \varphi^{\beta}}{\partial x ^{j}}=0, \quad \forall \alpha, \beta,
$$
which is the PHWC condition in coordinates, cf. \cite[Lemma 3]{lub}.

Now, due to the PHWC hypothesis, on $M$ we have a metric $f$-structure, $F$, such that $\varphi$ is $(F,J)$-holomorphic. Then for any $f$ holomorphic (locally defined on $N$), $f \circ \varphi$ is $f$-holomorphic. PHM condition requests that $f \circ \varphi$ be harmonic for all holomorphic function $f$. This means $(F\mathrm{div}F)(f\circ \varphi)=0, \forall f \text{holomorphic on} \ N$ (cf. Appendix), that is: $\mathrm{d}\varphi(F\mathrm{div}F)(f)=0, \forall f \ \text{holomorphic on} \ N$. But $N$ is a Hermitian manifold, so, according to Prop. A$_2$ cited in the Appendix, the above requirement is equivalent to $\mathrm{d} \varphi (F \mathrm{div} F)=0$ and, as $F\mathrm{div}F$ is horizontal, therefore equivalent to $F\mathrm{div}F=0$.

The proof of the converse is straightforward, taking into account that any $f$-holomorphic function on a manifold endowed with a cosymplectic $f$-structure is harmonic.
\end{proof}

Notice that the proof above shows us that in general a pseudo-harmonic morphism need not to be a harmonic map (but only PHWC).

\subsection{Foliations locally defined by pseudo-harmonic morphisms}
In the sequel we shall explore the conditions that a PHWC foliation must satisfy in order to be given locally by the fibers of a pseudo-harmonic morphism with values in $\CC^{n}$.

By analogy with the case of harmonic morphisms (see \cite {ud}), we introduce the following notion:

\begin{de}
Let $\FF$ be a codimension $2n$ foliation on a Riemannian manifold $(M, g)$. We say that $\FF$ {\it produces pseudo-harmonic morphisms} on $(M, g)$ if each point of $M$ has an open neighbourhood $U$ which is the domain of a submersive pseudo-harmonic morphism $\varphi: (U, g \vert _{U}) \longrightarrow (N^{2n}, J, h)$ whose fibers are open subsets of the leaves of $\FF$ (clearly, we may take $N^{2n}$ to be $\CC^{n}$).\end{de}

Note that the induced pseudo-harmonic morphism (PHM) is not unique, as the composition of a holomorphic map with a PHM is again a (local) PHM  (see also \cite[Prop. 3]{lub} for the strong PHM's case). 

\begin{pr}
On a Riemannian manifold $(M^{m}, g)$,
a PHWC foliation $\FF$ of codimension $2n$ $(m > 2n > 2)$ 
with integrable associated $f$-structure, $F$, produces (strong) pseudo - harmonic morphisms if and only if $F \mathrm{div} F=0$.
\end{pr}

\begin{proof}
According to Prop. 2.4, $\FF$ is locally defined by  $(F,J)$-holomorphic submersions onto some complex manifold $(N,J)$, which are in particular PHWC. Then, according to Prop. 3.2 above, these submersions are PHM's if and only if $F \mathrm{div} F=0$. 
\end{proof}

\medskip

The case of PHWC foliations of codimension 4 will be special, because, cf. \cite [p. 239]{ud}, any hermitian structure globally defined on $\RR^{4}$ is K\"ahler.

\medskip

Let us notice that the condition $F \mathrm{div} F=0$ is equivalent to: 
$$ 
2n(\mu^{\mathcal{H}})^{\flat}-(m-2n)(\mu^{\mathcal{V}})^{\flat}= \sum_{i=1}^{n}\mathrm{d}\Phi(e_{i}, Fe_{i}, \cdot),
$$ 
where $\Phi (X, Y) = g(X, FY)$ is the {\sf fundamental 2-form}.
\medskip

Recall now that, on a Riemannian  manifold $M^{m}$, a codimension $2n \neq 2$ conformal foliation $\FF$ (with the tangent distribution
$\VV$) produces {\it harmonic} morphisms if and only if the following 1-form is exact, cf. \cite{bry}:
\begin{equation} \label{harm}
W^{\flat} = (2n-2)(\mu^{\mathcal{H}})^{\flat} - 
(m - n)(\mu^{\mathcal{V}})^{\flat}
\end{equation}

In particular, let us notice that a transversally almost Hermitian foliation $\FF$, of codimension $2n$ $(m > 2n > 2)$ with associated $f$-structure, $F$, produces strong pseudo-harmonic morphisms if and only if $\mu^{\VV}=0$. But this conditions forces defining submersions to be moreover harmonic morphisms (the horizontal distribution is in this case totally geodesic, and the above criterion applies). 

\begin{re}
A foliation of even codimension, which produces harmonic morphisms, produces also "pure" {\it pseudo}-harmonic morphisms, by composing the distinguished submersions with (local) holomorphic, non-conformal maps.

Indeed, let $\pi: (M^{m}, g) \longrightarrow (N^{2n}, J, h)$ ($n \geq 2$) be a harmonic morphism onto a cosymplectic manifold (i.e. $\mathrm{div}J=0$). Then, for any holomorphic map $\varphi: (N, J, h) \longrightarrow (N^{\prime}, J^{\prime}, h^{\prime})$ to
a K\"ahler manifold, the composition 
$\widehat \varphi = \varphi \circ \pi$ is a {\it pseudo}-harmonic morphism. (This is a corollary of \cite[Prop. 3.4] {aab}.)
\end{re}

\subsection{The canonical foliation on almost contact manifolds}
In the rest of this section, we study on which (almost) contact metric manifolds, the foliation $\FF _{\xi}$ produces (pseudo -) harmonic morphisms (see also \cite {ian} for an independent approach). 
The notations are those one in \cite {ble}. Recall here only that 
$\Phi(X, Y) = g(X, \phi Y)$ is the \textit{fundamental 2-form} of an almost contact metric structure.

\begin{pr}
Let $(\phi, \xi, \eta, g)$ be an almost contact metric structure on a manifold $M^{2n+1}$.

\emph{(i)} If $(\phi, \xi, \eta, g)$ is {\bf nearly cosymplectic} 
(i.e. satisfies $(\nabla_{X} \phi)X=0, \forall X$), then $\FF _{\xi}$ produces harmonic morphisms with 1-dim. fibers, of Killing type (if $n > 1$). 

\emph{(ii)} If $(M^{2n+1}, \eta)$ is a {\bf K-contact manifold} (i.e. a contact manifold for which $\xi$ is Killing), the foliation associated to the Reeb vector field always produces harmonic morphisms with 1-dim. fibers, of Killing type ($n > 1$). 

\emph {(iii)} If $\phi$ is {\bf normal} and $\delta \Phi = 0$ (equivalently: $\phi \mathrm{div}^{\DD} \phi =0$),
then $\FF _{\xi}$ produces \emph {pseudo}-harmonic morphisms.

\emph{(iv)} If $(\phi, \xi, \eta, g)$ is
{\bf quasi - Sasakian} (i.e. is normal and satisfies $\mathrm{d}\Phi = 0$), then the foliation
$\mathcal{F}_{\xi}$ produces harmonic morphisms with 1-dim. fibers, of Killing type (if $n > 1$). 

\emph{(v)} If $(\phi, \xi, \eta, g)$ is {\bf $\alpha$ - Sasakian} (i.e. satisfies  
$\left(\nabla_{X}\phi \right)Y = \alpha [g(X, Y)\xi - \eta(Y) X]$, $\forall X, Y$,
for a function $\alpha$ on M),
then the foliation $\FF_{\xi}$ produces harmonic morphisms with 1-dim. fibers, of Killing type. 

\emph{(vi)} If $(\phi, \xi, \eta, g)$ is {\bf Kenmotsu} (i.e. satisfies 
$\left(\nabla_{X}\phi \right)Y = g(\phi X, Y)\xi - \eta(Y)\phi X$, $\forall X, Y$),
then the foliation $\FF_{\xi}$ produces harmonic morphisms with 1-dim. fibers, of warped product type. 

\end{pr}

\begin{proof}
In this context, $\HH= \DD =\Ker \eta$ and $\VV = Sp\{ \xi \}$.

$(i)$ On a nearly cosymplectic manifold, $\xi$ is automatically a Killing vector field, so $\FF _{\xi}$ is Riemannian, which is equivalent to $\mathcal{D}$ being totally geodesic. Also, from $(\nabla_{\xi} \phi)\xi=0$, we get $\nabla_{\xi}\xi=0$ so $\mu^{\VV}=0$. We conclude that $W=0$ and the assertion {\it (ii)} is obvious.

$(ii)$ In this case, $\xi$ being Killing, $\FF _{\xi}$ is again a Riemannian foliation (so $\DD$ is totally geodesic) and $\mu^{\mathcal{V}}=0$ (the integral curves of $\xi$ are
even geodesics, because on any contact manifold one has $\mathcal{L}_{\xi}\eta=0$). 
\footnote{Note that, even if for a contact manifold one can prove $\mu^{\DD}=0$ and $\mu^{\VV}=0$ 
(so $W=0$), $\FF _{\xi}$ is not conformal in general, so does not produce harmonic morphisms. It produces neither pseudo-harmonic morphisms, because $\phi$ is not generally a basic tensor field.}

$(iii)$ On a normal almost contact manifold, $\phi$ is integrable as mentioned above (Example 2.1). The normality of $\phi$ implies $\mathcal{L}_{\xi}\eta=0$ which forces the leaves of $\FF _{\xi}$ to be minimal (see Remark 4.1 bellow). Therefore the hypothesis $\delta \Phi = 0$ suffices to assure $F \mathrm{div} F=0$. 

$(iv)$ It is a standard fact that on a quasi - Sasakian manifold $\xi$ is again a Killing vector field (so, as above $\FF _{\xi}$ is Riemannian and $\DD$ is totally geodesic). The normality $\phi$ forces the leaves of $\FF _{\xi}$ being minimal (again from $\mathcal{L}_{\xi}\eta=0$). As above one see that $W=0$.

$(v)$ The relation from the definition of $\alpha$-Sasakian
structures (which implies $\phi$ normal!) gives us immediately that: $\eta (\nabla_{X}Y)=-\alpha g(X, \phi Y)$,
$\forall X, Y \in \Gamma(\DD)$. Therefore 
$B^{\DD}(X,Y)=2^{-1}\eta (\nabla_{X}Y + \nabla_{Y}X)\xi =0$. So $\DD$ is totally geodesic which is equivalent to $\FF _{\xi}$ being a Riemannian foliation.

Also, it is easy to see from definition that $(\nabla_{\xi}\xi)^{\DD}= 0$, so $\mu^{\VV}=0$ and then $W=0$.

$(vi)$ On a Kenmotsu manifold, $\nabla_{\xi}\xi = 0$ (so $\mu^{\VV}=0$) and $\mathcal{D}$ is an umbilic (integrable) distribution (with $\mu^{\DD} = -\xi$), so $\FF_{\xi}$ is conformal. Because, in this case, $d \eta = 0$, $\eta$ is locally exact, so $W$ is (locally) a gradient. Moreover the local warped-product structure of any Kenmotsu manifold gives us the type of the induced harmonic morphisms with 1-dim. fibers.
\end{proof}

\begin{co}
Let $(\phi, \xi, \eta, g)$ be a normal {\bf semi-cosymplectic} almost contact metric structure,
so in the class $\mathcal{C}_{3} \oplus \mathcal{C}_{7} \oplus \mathcal{C}_{8}$, according 
to \emph{Chinea -Domingo classification}, \cite {chin}.
Then $\FF _{\xi}$ produces pseudo-harmonic morphisms with minimal horizontal 
distribution. 

\end{co}

\begin{proof} 
This class of normal almost contact structures is defined by the following additional properties: $\delta \Phi =0$ and $\delta \eta= 0$. So ($iii$) above applies.
Note that the condition $\delta \eta= 0$ is equivalent to $\DD$ minimal.
\end{proof} 

\subsection{An example of pseudo-harmonic morphism that
has $(1,2)$-symplectic associated $f$-structure, but it is not PHH} The following considerations argue that \cite [Prop. 7]{mo} is really a generalization of \cite[Theorem 4.1] {aab} (which was not justified by an example so far).

\begin{lm}
On a K-contact manifold which is not Sasaki, $\phi$ is
an $(1,2)$-symplectic $f$-structure that is not parallel in horizontal (i.e. tangent to $\DD$) directions,
with respect to $\nabla^{\mathcal{D}}$.
\end{lm}

\begin{proof}
The fact that $\phi$ is an $(1,2)$-symplectic $f$-structure 
(i.e. $\left(\nabla_{X}\phi \right)Y+\left(\nabla_{\phi X}\phi \right)\phi Y$ $=0$) 
follows immediately from
Olszak's formula, \cite{ol}, which holds on any contact manifold:
\begin{equation*}
\left(\nabla_{X}\phi \right)Y+\left(\nabla_{\phi X}\phi \right)\phi Y=2g(X, Y)\xi-\eta(Y)(X+hX+\eta(X)\xi).
\end{equation*}

For the second statement we shall use the following relation, which holds on any contact manifold too
(cf. \cite [Lemma 6.1]{ble}):
$$
2g\left((\nabla_{X}\phi)Y, Z \right)=g(N^{(1)}(Y, Z), \phi X)+
2d\eta(\phi Y, X)\eta(Z)-2d\eta(\phi Z, X)\eta(Y).
$$

Suppose that $\nabla^{\mathcal{D}}_{X} \phi=0,
\forall X \in \Gamma(\mathcal{D})$.
Then, from the above relation, for $X \in \Gamma(\DD)$ arbitrary, we derive: 
$g(N^{(1)}(Y, Z), \phi X)=0, \forall Y, Z \in \Gamma(\mathcal{D})$. In turn, this implies:
 $\left( N^{(1)}(Y, Z) \right)^{\mathcal{D}}=0, \forall Y,Z \in \Gamma(\DD)$.

It is easy to verify that also the component collinear with $\xi$ of $N^{(1)}(Y, Z)$ is zero:
 $g(N^{(1)}(Y, Z), \xi) = \eta([\phi X, \phi Y]-[X, Y]) = 0$ 
(because we are on a contact manifold and $N^{(2)}=(\mathcal{L}_{\phi X}\eta)(Y)-
(\mathcal{L}_{\phi Y}\eta)(X)=0$).

Now, $N^{(1)}(Y, \xi)=\phi\left((\mathcal{L}_{\xi} \phi)X\right)$ and the K-contact hypothesis 
assures that this component of $N^{(1)}$ is zero, too. So, in our assumption, we have obtained $N^{(1)}=0$, 
which can not be true, because our manifold is not Sasaki.
\end{proof}

\begin{co}
On a K-contact manifold which is not Sasaki, $\FF _{\xi}$ produces a pseudo-harmonic morphism that has $(1,2)$-symplectic associated $f$-structure, but it is not PHH.
\end{co}

\begin{proof}
We have to apply the above Lemma, together with Remark 3.1 and Proposition 3.2. 
When applying Remark 3.1, let $J$ be the projection of $\phi$ which is, in particular, cosymplectic as we can easily prove.
\end{proof}

\section{Pseudo-harmonic morphisms with one-dimensional fibers}

Let $\varphi: (M^{2n+1},g) \longrightarrow (N^{2n}, J, h)$ be a PHWC submersion with one dimensional fibers onto a (almost) Hermitian manifold. 
We have seen that the horizontal bundle $\HH$ inherits an almost complex structure $J_{\HH}$, so in particular is oriented.
If $M$ is orientable, an assumption which we adopt from now on,
then the quotient bundle $\VV=TM / \HH$ is an orientable line bundle over $M$, hence trivial.
Therefore it admits a globally defined nowhere vanishing section
(when $M$ is oriented and compact, the existence of this field forces the Euler characteristic of $M$ to be zero). Choosing a unitary section $\xi$ of this type, one obtain an almost contact metric structure $(\phi, \xi, \eta, g)$ on $M$. 
Indeed, if we take $\phi | _{\HH}=J_{\HH}$, $\phi | _{\VV}=0$ and
$\eta(X)=g(X, \xi), \forall X \in \Gamma(TM)$, then all the conditions in the definition of an almost contact metric structure are satisfied
(i.e. $\phi^{2}=-I+\eta \otimes\xi$, $\eta (\xi) =1$
and $g(\phi X, \phi Y)=g(X, Y)-\eta(X)\eta(Y)$,
cf. \cite {ble}). Moreover, $\varphi$ is then a $(\phi, J)$-holomorphic map (i.e. $d\varphi \circ \phi = J \circ d\varphi$).

\medskip

In the following we shall answer de question: under which additional hypothesis, $\phi$ is normal? 

\medskip

Let us now recall some standard facts on 1-dimensional foliations, pointed out in \cite [Section 10.5]{ud}.

\begin{re}
The 1-dimensional foliation $\mathcal{F}_{\xi}$ has minimal (i.e. geodesic) leaves (with respect to $g$) if and only if $\mathcal{L}_{\xi}\eta=0$ (which is equivalent to 
$\mathcal{L}_{\xi}\Gamma(\mathcal{H}) \subset \Gamma(\mathcal{H})$ as well as to $d\eta(\xi, X)=0$, for any horizontal field $X$).

The \emph{integrability tensor of the horizontal distribution} is: $I^{\mathcal{H}}(E, F) = [E^{\mathcal{H}}, F^{\mathcal{H}}]^{\mathcal{V}}$, $\forall E, F \in \Gamma(TM)$. Then $\Omega :=d\eta$ will be called \emph{integrability 2-form of} $\mathcal{F}_{\xi}$. For $X, Y$ horizontal, we have:
$$
\Omega(X, Y)=-\eta([X, Y])=-g(I^{\mathcal{H}}(X, Y), \xi).
$$
So $I^{\mathcal{H}}(X, Y)=-\Omega(X, Y)\xi$. In addition, for a foliation with minimal leaves, $\Omega$ is a \emph{basic} form (i.e. $\imath _{\xi}\Omega=0$ and
$\mathcal{L}_{\xi}\Omega=0$).

Note also that the natural decomposition of the metric: $g=g^{\mathcal{H}}+g^{\mathcal{V}}$ becomes in this case:
$g=g^{\mathcal{H}}+\eta \otimes \eta$. 
\end{re}

\begin{te}
Let $\varphi: (M^{2n+1},g) \longrightarrow (N^{2n}, J, h)$ be a PHWC submersion from an oriented Riemannian manifold onto a Hermitian manifold.

The almost contact structure naturally inherited by M is normal if and only if (one of) the following equivalent conditions hold:

\medskip
\emph {(i)} the fibers are minimal and $d\eta (\phi X, \phi Y) = d\eta (X, Y), \forall X, Y$ $($the integrability form $\mathrm{d}\eta$ is $\phi$-invariant, or it is a $2$-form of bidegree $(1, 1))$

\medskip
\emph {(ii)} $\nabla_{\xi} \phi = 0$ $(\phi$ is parallel along the fibers$)$. In particular, the fibers are minimal.

\end{te}

\begin{proof}
The induced almost contact structure ($\phi, \xi, \eta$) is normal if and only if
$$
N^{(1)}(X, Y)=[\phi, \phi](X, Y)+2d\eta(X, Y)\xi=0.
$$
For a couple of vertical and horizontal vector fields (where $\HH = \DD = \Ker \eta$) we have:
$$
N^{(1)}(X, \xi)=[\phi, \phi](X, \xi)+2d\eta(X, \xi)\xi=
\phi \left( (\mathcal{L}_{\xi} \phi)X \right)-(\mathcal{L}_{\xi} \eta)(X) \xi.
$$ 
 
We have already noted in Lemma 2.1 that 
$(\mathcal{L}_{\xi} \phi)(X) ^{\HH}=0, \forall X \in \Gamma(\HH)$. 
So $N^{(1)}(X, \xi)=0, \forall X \in \Gamma(\HH)$ is equivalent with the minimality of fibers.

On the other hand, for a couple of horizontal vector fields we have:
$$
N^{(1)}(X, Y)=[\phi, \phi](X, Y)^{\HH} + 
\left( \eta([\phi X, \phi Y]) - \eta([X, Y]) \right) \xi.
$$ 
But, according to Lemma 2.2, we have already:
$$
[\phi, \phi](X, Y)^{\mathcal{H}}=0, \forall X, Y \in \Gamma(\HH).
$$

So $N^{(1)}(X, Y)=0, \forall X, Y \in \Gamma(\HH)$ is equivalent to: 
\begin{equation} \label {finv}
\eta([\phi X, \phi Y]) = \eta([X, Y]), \forall X, Y \in \Gamma(\HH),
\end{equation}
which is equivalent to $d\eta (\phi X, \phi Y) = d\eta (X, Y), \forall X, Y$. We have proved that the condition $\phi$ normal is equivalent to ($i$).

Note that \eqref {finv} is equivalent also with
$N^{(2)}(X, Y) = 0, \forall X, Y \in \Gamma(\HH)$. 

In order to finalize the proof, we have to "calculate" $\nabla_{\xi}\phi$.
According to \cite [Lemma 6.1]{ble}, one has:
$$
2g\left( (\nabla_{\xi}\phi)X, Y \right)= 3d\Phi(\xi, \phi X, \phi Y)-
3d\Phi(\xi, X, Y)+ N^{(2)}(X, Y), \forall X, Y \in \Gamma(\HH),
$$
where, as usual, $\Phi(X, Y):=g(X, \phi Y)$.

Now, precisely because $\phi$ is a basic tensor, that is
 $(\mathcal{L}_{\xi} \phi)(X)^{\HH}=0, \forall X \in \Gamma(\HH)$,
we can check directly (or by applying Corollary 2.1,(i)) that:
$$
d\Phi(\xi, \phi X, \phi Y)-d\Phi(\xi, X, Y)=0.
$$

Therefore
$N^{(2)}(X, Y)=0, \forall X, Y \in \Gamma(\HH)$
is equivalent to 
$\left((\nabla_{\xi}\phi)X\right)^{\mathcal{H}}=0$.

It is easy to prove that minimality of the fibers assures both $\left((\nabla_{\xi}\phi)X\right)^{\mathcal{V}}=0, \forall X \in \Gamma(\HH)$ and  $(\nabla_{\xi}\phi)(\xi)=0$. 

Using ($i$), we can conclude that the normality of $\phi$ is equivalent also with $\nabla_{\xi}\phi=0$, that is ($ii$).
\end{proof}

Harmonicity of a PHWC submersive map does not interfere, {\it a priori}, with the normality of $\phi$ (harmonicity reduces to a condition in the horizontal bundle, $\phi (\mathrm{div}^{\HH}\phi)=0$, while normality is coming from the parallelism of $\phi$ along the fibers).

\begin{co}
Let $\varphi: (M^{2n+1},g) \longrightarrow (N^{2n}, J, h)$ be a PHWC submersion from an oriented Riemannian manifold, into a Hermitian manifold. If the fibers are all diffeomorphic to circles and the following relation is satisfied:
\begin{equation}\label{CR}
[\dif \varphi \circ \nabla^{M}_{\xi} \circ \dif \varphi^{*}, J]=0,
\end{equation}
then the induced almost contact structure on $M$ is normal and $\varphi$ is a \emph{contact bundle}, cf. \cite{mor}, \cite{mori}.
\end{co}

\begin{proof}
Notice that \eqref{CR} translates $\nabla^{M}_{\xi}\phi=0$ on $\HH$ and the above theorem assures us that the almost contact structure induced on $M$ is normal. But \cite[Corollary 1]{mori}, shows us that {\it any} normal almost contact manifold  $(M, \phi, \xi, \eta)$ such that  $\xi$ is a {\it closed} vector field is a principal circle bundle over a complex manifold with $\eta$ a connection form and $\xi$ a vertical vector field corresponding to the unit vector in $L(S^{1})$, which is called a contact bundle (instead of $\xi$ closed, one can suppose $\xi$ regular and $M$ compact).
\end{proof}

\section{Pseudo-harmonic morphisms with two-dimensional fibers}

First, we shall remark that any PHWC map with 1-dimensional fibers from an orientable Riemannian manifold to a Hermitian manifold gives rise to a naturally associated PHWC map with 2-dimensional fibers.

We recall here that \emph{the cone} $\mathcal{C}(M)$ over an almost contact metric manifold $(M^{2n+1}, \phi, \xi, \eta, g)$ is simply $M^{2n+1}\times \mathbb{R}$ with an almost complex structure defined by:
$$
J\left(X, f \frac{d}{dt} \right)=(\phi X-f\xi, \eta(X)\frac{d}{dt}).
$$
We endow the cone with a warped product metric $\widehat{g}=dt \otimes dt + t^{2}g$. Let $\pi$
denote de standard projection $\pi: \mathcal{C}(M)\longrightarrow M$.

\begin{pr}
Let $\varphi: (M^{2n+1},\phi, \xi, \eta, g) \longrightarrow (N^{2n}, J^{N}, h)$ $(n \geq 2)$ be a smooth map from an almost contact metric manifold to an (almost) Hermitian manifold. Denote by $\widehat{\varphi}$ the composition map $\varphi \circ \pi$.

\emph{(i)} $\varphi$ is $(\phi, J^{N})$-holomorphic if and only if $\widehat{\varphi}$ is holomorphic. 

\emph{(ii)} $\varphi$ is a PHWC harmonic map (strong PHM) if and only if so is $\widehat{\varphi}$ \ (if $h$ is moreover a K\"ahler metric on $(N, J^{N})$).

\emph{(iii)} $\varphi$ is submersive if and only if so is $\widehat{\varphi}$. In this case, if the inverse image by $\varphi$ of a complex submanifold is a minimal one, the same is true for $\widehat{\varphi}$.

\emph{(iv)} $\varphi$ is PHM (pulls holomorphic functions back to harmonic ones) if and only if so is $\widehat{\varphi}$.
\end{pr}

\begin{proof}
$(i)$ Suppose that $\varphi$ is $(\phi, J^{N})$-holomorphic. Then, we have:
$$
d\widehat{\varphi}\left(J\left(X, f \frac{d}{dt} \right)\right)=d\varphi(\phi X -f\xi)
= d\varphi(\phi X) = 
$$
$$
J^{N}d\varphi(X)=J^{N}\left(d\widehat{\varphi}\left(X, f \frac{d}{dt} \right)\right).
$$
Conversely, if $\widehat{\varphi}$ is holomorphic then we have:
$d\varphi(\phi X) - f d\varphi(\xi)=J^{N}d\varphi(X)$. But $d\varphi(\xi) = -J^{N}d\widehat{\varphi}\left(0, \frac{d}{dt} \right) = 0$ and our assertion follows.

$(ii)$ We have only to remark that $\pi$ is a horizontally homothetic harmonic morphism. This is because it is obviously horizontally homothetic and we can easily check that its fibers are geodesics, so minimal. 

So we are in the hypothesis of \cite[Prop. 3.4] {aab},  which states 
that $\varphi$ is a PHWC harmonic map if and only if so is $\widehat{\varphi}$.

$(iii)$ Because $\pi$ is submersive (it is a harmonic morphism with the dimension of the target space greater than 4), the first assertion is elementary. Using the fact that $\pi$ is horizontally homothetic, one can verify directly the second assertion. 

$(iv)$ We can directly check the equivalence using that $\pi$ is harmonic morphism. 
\end{proof}

\bigskip

Now, let $\varphi: (M^{2n+2},g) \longrightarrow (N^{2n}, J, h)$ be a PHWC submersion with 2 - dimensional fibers, from an oriented Riemannian manifold into a K\"ahler manifold. Then, on $M$, we shall have two almost Hermitian structures with respect to which $\varphi$ is holomorphic (this will be called \emph{adapted almost Hermitian structure} as in \cite {ud}).
Indeed, consider the natural split induced by the submersion 
$TM = \mathcal{H} \oplus \mathcal{V}$. On $\mathcal{H}$ we have the almost complex 
structure $J_{\mathcal{H}}$ (and therefore an orientation),
 compatible with the metric, induced by the PHWC condition. 
On $\mathcal{V}$, consider the orientation such that
$\mathcal{H} \oplus \mathcal{V}$ has the orientation of $M$ and also let $J_{\mathcal{V}}$ 
be the almost Hermitian structure given by 
rotation through $+ \pi/2$. Then we have two \emph{adapted} almost Hermitian structure on $M$:
$J^{+}=(J_{\mathcal{H}}, J_{\mathcal{V}})$ and $J^{-}=(J_{\mathcal{H}}, -J_{\mathcal{V}})$.

The first fact to be remarked is the following (which appear, in other context, also in \cite {ud}):

\begin{re} 
Under the above hypothesis, we have:
$$
\nabla_{E}^{\mathcal{V}}J_{\mathcal{V}}=0, \quad \forall E \in \Gamma(TM),
$$
where $\left(\nabla_{E}^{\mathcal{V}}J_{\mathcal{V}}\right)(V):=
\left(\nabla_{E} J_{\mathcal{V}} V\right)^{\mathcal{V}}
-J_{\mathcal{V}}\left(\nabla_{E} V \right)^{\mathcal{V}}, 
\forall V \in \Gamma(\mathcal{V})$.
\end{re}

\begin{proof}
As in Proposition 2.1, we can prove that:

$g \left( (\nabla_{E}J_{\mathcal{V}})(V), V \right)=
g\left( (\nabla_{E}J_{\mathcal{V}})(V), J_{\mathcal{V}}V \right)=0$.

Because, in this case, $V$ and $J_{\mathcal{V}}V$ form a basis in $\mathcal{V}$, the conclusion follows. 
\end{proof}

Now, we assert that \cite [Theorem 5.1] {gud} holds without the 
horizontally conformal hypothesis and we state it in our context. Recall that an almost complex submanifold $P$ of an almost Hermitian manifold $(M, J, g)$ is \emph{superminimal} if $J$ is parallel
along $P$ i.e. $\nabla_{V}J=0$ for all vector fields $V$ tangent to $P$.

\begin{te}
Let $\varphi: (M^{2n+2},g) \longrightarrow (N^{2n}, J, h)$ be a PHWC submersion with complex 1-dimensional fibers, from an oriented Riemannian manifold into a Hermitian manifold. Consider on $M$ one of the induced almost complex structures $(J_{\HH},\pm J_{\VV})$, compatible with $g$. If  the fibers are superminimal, then this almost complex structure is integrable. 
\end{te}

\begin{proof}
Consider on $M$ the almost Hermitian structure $(J^{+}, g)$ (the other case is completely analogous). Let $N^{J^{+}}$ be its Nijenhuis torsion.

Recall that Lemma 2.2 assures us that $N^{J^{+}}(X,Y)^{\HH}=0, \forall X,Y \in \Gamma(\HH)$. We know also that $J_{\HH}$ is basic, so: $[(\mathcal{L}_V J^{+})(X)]^{\HH}=0, \forall V \in \Gamma(\VV), X \in \Gamma(\HH)$. Then it is easy to verify thet $J^{+}$ integrable (i.e. $N^{J^{+}}=0$) will be equivalent to the following conditions:
\begin{equation*}
\begin{split}
&(\alpha) \ N^{J^{+}}(X,Y)^{\VV}=0, \qquad
(\beta) \ N^{J^{+}}(X, V)^{\VV}=0, \quad \forall X,Y \in \Gamma(\HH), V \in \Gamma(\VV).
\end{split}
\end{equation*}

According to \cite [Lemma 3.1]{sve}, on the almost Hermitian manifold $(M, J^{+}, g)$, for any section $X$ of $\mathcal{H}$ and any vector $V$ tangent to $\mathcal{V}$, we have:
\begin{equation*}
g(\nabla_{J^{+} Y}X+\nabla_{Y}J^{+} X, V)=g(X, (\mathcal{L}_{V}J^{+}-\nabla_{V}J^{+})Y), 
\quad \forall Y \in \Gamma(TM).
\end{equation*}

Taking into account that $J_{\HH}$ is basic and $[(\nabla_{V} J^{+})Y]^{\HH}=0$ (from the superminimality hypothesis), the above formula, for any $X, Y \in \Gamma (\HH)$, reduces to:
\begin{equation*}
g(\nabla_{J^{+} Y}X + \nabla_{Y}J^{+} X, V)=0, \forall V \in\Gamma(\VV).
\end{equation*}

In particular we have also: $g([X, J^{+} Y]+[J^{+} X, Y], V)=0.$ From this we get immediately $(\alpha)$.

Now, let us notice that the following formula holds:
\begin{equation*}
g\left(N^{J^{+}}(X, V), V \right)=
-g\left(J^{+}(\nabla_{V}J^{+})(X)-(\nabla_{J^{+}V}J^{+})(X), V\right),
\end{equation*}
and a similar one for $g(N^{J^{+}}(X, V), J^{+}V)$, where $\{V, J^{+}V \}$ is now an orthonormal frame in $\VV$. As we suppose $[(\nabla_{V} J^{+})X]^{\VV}=[(\nabla_{J^{+}V} J^{+})X]^{\VV}=0$, it is clear that $(\beta)$ is satisfied too.

\end{proof}

Remark that, contrary to the case treated by Theorem 4.1, in general the integrability of the induced almost complex structure does not imply (super)minimality of the fibers (to put it otherwise, superminimality hypothesis is too strong).

\section{Appendix}

A function $f: (M,J) \rightarrow \CC$ on an almost complex manifold is called holomorphic iff $\dif f \circ J = \mathrm{i} \cdot \dif f$.

It is a classical fact that:
\begin{quote}
$\textbf{A}_1 .$ \textit{A holomorphic function on a cosymplectic (almost Hermitian) manifold is always harmonic.}
\end{quote}
\begin{proof}
Let $(M,g, J)$ an almost Hermitian cosymplectic manifold and $f: M \rightarrow \CC$ a holomorphic function. In this case the Hodge Laplacian is identical to the rough Laplacian. Holomorphicity of $f$ translates as follows:
$\dif f (JX) = \mathrm{i}\dif f(X)$ or $JX(f)=\mathrm{i}X(f)$. Letting $\{E_i, JE_i\}_{i=\ov{1,n}}$ an adapted orthonormal frame, we shall have:

\begin{equation*}
\begin{split}
\tr \nabla ^2 (f)&=\sum_{j} E_j(E_j(f))+ JE_j(JE_j(f))- \left(\nabla_{E_j} E_j + \nabla_{JE_j} JE_j \right)(f)\\
&=\sum_{j} -E_j(J^2 E_j(f))+ \mathrm{i}JE_j(E_j(f))+ J^2 \left(\nabla_{E_j} E_j + \nabla_{JE_j} JE_j \right)(f)\\
&=\sum_{j} -\mathrm{i}E_j(J E_j(f))+ \mathrm{i}JE_j(E_j(f))+ \mathrm{i}J \left(\nabla_{E_j} E_j + \nabla_{JE_j} JE_j \right)(f)\\
&=\sum_{j} -\mathrm{i}\left[ [E_j,J E_j] - J \left(\nabla_{E_j} E_j + \nabla_{JE_j} JE_j \right)\right](f)\\
&=-\mathrm{i}(\mathrm{div}J)(f)=-(J\mathrm{div}J)(f)= \ 0.
\end{split}
\end{equation*}
\end{proof}

Of course this is a particular case of the result of {\it Andr\'e
Lichnerowicz} ('70), \cite{lih}: 
a (anti)holomorphic map from a {\it cosymplectic (semi-K\"ahler)} manifold to a (1, 2)-{\it symplectic (quasi-K\"ahler)} manifold is harmonic.

If $J$ is integrable we have a converse result: 

\begin{quote}
$\textbf{A}_2 .$ \textit{A Hermitian manifold $M$ is cosymplectic if and only if every holomorphic function on an open subset of $M$ is harmonic}.
\end{quote}

The proof of this second statement is based upon the Newlander - Nirenberg theorem which assures us that if $J$ is integrable, then each point of $M$ has a neighbourhood in which there exists $\mathrm{dim}_{\RR}M$ (functionally) independent holomorphic functions (so cancelation of a vector field on any holomorphic function suffices to state that the vector field is zero).

\medskip

Analogously, a function $f: (M,F) \rightarrow \CC$ on a manifold endowed with a $f$-structure is called \textit{f-holomorphic} iff $\dif f \circ F = \mathrm{i} \cdot \dif f$. In particular, $\dif f (V)=0, \forall V \in \Ker F$.

Extending the discussion presented in \cite{slob} for normal almost contact structures, we can prove the analogous version of the above two results, for $f$-holomorphic functions on manifolds endowed with an integrable / cosymplectic $f$-structure (recall that a $f$-structure is called \textit{cosymplectic} iff $F\mathrm{div}F=0$).

\section {Gloss}

In 1848, Jacobi \cite{jac} studied the following problem for $m = 3$:
\begin{quote}
\textit{Let $\phi : U \rightarrow \CC$ be a smooth map from an open subset of $\RR^m$. Under what conditions on $\phi$ is the composition $f \circ \phi$ harmonic for an arbitrary holomorphic function $f : V \rightarrow \CC$ defined on an open subset of $\CC$}?
\end{quote}

Now the answer is well known: since any harmonic function is locally the real part of a holomorphic function it is clear that $\phi$ satisfies Jacobi's condition if and only if it is a harmonic morphism. This was one of starting points for the history of the long debated notion of harmonic morphism.

We can see the present discussion on \textit{pseudo}-harmonic morphisms as the study, in the most general setting, of \textit{the same} question. If $\phi$ is a smooth map from a Riemannian to a Hermitian manifold, then the answer is: $\phi$ must be PHWC map with cosymplectic induced $f$-structure. Then let us remark that, in fact, this approach (\textit{pulling back holomorphic functions defined on a complex manifold}
\footnote{The holomorphic functions on the codomain are not harmonic in general, as it was the case in the original context or in the strong PHM context.}) does not provide purely harmonic functions, but functions which are ($f$-)holomorphic too and, moreover, harmonic due to the geometry of the domain.

We have considered mainly the case when the almost complex structure on the codomain is integrable, as in the nonintegrable case the existence of holomorphic functions becomes problematical.

\bigskip
{\bf Acknowledgements}.
I am grateful to Professor Vasile Brinzanescu, for his special generous patience to pursue this work and advise me. A significant step to the final version of this work has been done in the frame of my PostDoc Program at {\it Laboratoire de Math\'ematique de Brest, CNRS - UMR 6205}, France. In this respect, I am indebted to Professor Eric Loubeau for his warm collaboration.

For the present revisited version, I am grateful to Radu Pantilie for enlightening discussions and I acknowledge the partial support by the CEx grant no.\ 2-CEx 06-11-22/25.07.2006.

\end{document}